\documentclass[11pt]{article}
\usepackage{amsmath,latexsym,amssymb}
\usepackage[draft]{hyperref}
\def\leqst{\leq_{\mathrm{st}}}
\def\geqst{\geq_{\mathrm{st}}}
\def\leqlr{\leq_{\mathrm{lr}}}

\def\hmlr{\sim_{\mathrm{hmlr}}}
\def\given{\hspace{0.8pt}|\hspace{0.8pt}}

\newcommand{\mfalls}{\quad\mbox{if \;}}
\newcommand{\msonst}{\quad\mbox{otherwise}}
\renewcommand{\cases}[1]{\left\{\begin{array}{rl}#1\end{array}\right.}

\def\l{\ell}
\newcommand{\mbu}{\quad\mbox{ and }\quad}

\newcommand{\mbs}[1]{\mbox{ \;#1\; }}
\newcommand{\mbsl}[1]{\mbox{ \;#1}}

\newcommand{\mf}{\quad\mbox{\;for \;}}
\newcommand{\mfa}{\quad\mbox{\;for all \;}}

\newcommand{\mfasts}{\quad\mbox{almost surely}}

\newcommand{\mfs}{\quad\mbox{a.s.}}

\renewcommand{\P}{\mathbf{P}}
\newcommand{\E}{\mathbf{E}}

\newcommand{\R}{{\mathbb{R}}}
\newcommand{\N}{{\mathbb{N}}}
\newcommand{\CB}{{\mathcal{B}}}
\newcommand{\CC}{{\mathcal{C}}}

\newcommand{\ve}{\varepsilon}
\def\hyp{\mathrm{hyp}}
\newcommand{\Poi}{\mathrm{Poi}}
\def\d{\mathrm{d}}

\newcommand{\equ}[1]{(\ref{#1})}
\providecommand{\qed}{{}\hfill {}\hfill{$\Box
$}\vspace{0.3cm}\pagebreak[2]\par}
\newcommand{\Section}[1]{\section{#1}\setcounter{equation}{0}}

\newcommand{\DS}{\displaystyle}

\newtheorem{theorem}{Theorem}
\newtheorem{proposition}{Proposition}[section]
\newtheorem{lemma}[proposition]{Lemma}
\newtheorem{uremark}[proposition]{Remark}
\newtheorem{definition}[proposition]{Definition}

\newenvironment{remark}{\begin{uremark}\rm}{\end{uremark}}

\def\CH{\mathcal{H}}

\def\CB{\mathcal{B}}

\def\ARG{\;\boldsymbol{\cdot}\;}
\begin{document}
\title{Stochastic ordering of classical discrete distributions}
\author{Achim Klenke\\
Johannes Gutenberg-Universit{\"a}t Mainz\\
Institut f{\"u}r Mathematik\\
Staudingerweg 9\\
55099 Mainz\\
Germany\\
\small
math@aklenke.de\\
\and
Lutz Mattner\\
Universit{\"a}t Trier\\
FB IV - Mathematik\\
54286 Trier\\
Germany\\
mattner@uni-trier.de
}

\date{\tiny
Submitted 07 March, 2009\\Revised 07 March, 2010}
\maketitle
\begin{abstract}
For several pairs $(P,Q)$ of classical distributions on $\N_0$, we show that their stochastic ordering $P\leqst Q$ can be characterized by their extreme tail ordering equivalent to $ P( \{k_\ast \})/Q(\{k_\ast\}) \geq 1 \geq  \lim_{k\rightarrow k^\ast} P(\{k\})/Q(\{k\})$, with $k_\ast$ and $k^\ast$ denoting the minimum and the supremum of the support of $P+Q$, and with the limit to be read as $P(\{k^\ast\})/Q(\{k^\ast\})$ for $k^\ast$ finite. This includes in particular all pairs where $P$ and $Q$ are both binomial ($b_{n_1,p_1} \leqst b_{n_2,p_2}$ if and only if  $n_1\le n_2$ and $(1-p_1)^{n_1}\ge(1-p_2)^{n_2}$,  or $p_1=0$), both negative binomial ($b^-_{r_1,p_1}\leqst b^-_{r_2,p_2}$ if and only if $p_1\geq p_2$ and $p_1^{r_1}\geq p_2^{r_2}$), or both hypergeometric with the same sample size parameter. The binomial case is contained in a known result about Bernoulli convolutions, the other two cases appear to be new.

The emphasis of this paper is on providing a variety of different methods of proofs: (i) half monotone likelihood ratios, (ii) explicit coupling, (iii) Markov chain comparison, (iv) analytic calculation, and (v) comparison of L{\'e}vy measures. We give four proofs in the binomial case (methods (i)-(iv)) and three in the negative binomial case (methods (i), (iv) and (v)). The statement for hypergeometric distributions is proved via method (i).
\end{abstract}
\vfill
\footnoterule\smallskip
\Section{Introduction}
\label{S1}
\subsection{Stochastic Ordering}
\label{S1.1}
For probability measures $P$ and $Q$ on the real numbers, the
stochastic ordering is the partial ordering
$$P\leqst Q \quad\Longleftrightarrow\quad
P([x,\infty))\leq Q([x,\infty)) \mbs{for all}x\in\R.
$$
This condition is equivalent to the existence of two real-valued
random variables $X$ and $Y$ with distributions $P$ and $Q$,
respectively, and such that $X\leq Y$ almost surely. In fact, let
$F^{}_P$ and $F^{}_Q$ denote the distribution functions of $P$ and
$Q$, respectively, and let $F_P^{-1}$ and $F_Q^{-1}$ be their
left-continuous inverses. That is,
$$
F_P^{-1}(t):=\inf\{x\in\R:\,F^{}_P(x)\geq t\}.
$$
Further, let $U$ be uniformly
distributed on $(0,1)$. Then $X:=F_{P}^{-1}(U)$ and
$Y:=F_Q^{-1}(U)$ have the desired property. Such a pair $(X,Y)$
is called a coupling.

Recall that $P\leqst Q$ is equivalent to the condition
that for any bounded and monotone increasing function $f:\R\to\R$,
we have
$$\int f\,\d P\;\leq\;\int f\,\d Q.$$
If $P$ and $Q$ have finite expectations, then taking
$f(x)=\max(-n,\min(x,n))$ and letting $n\to\infty$ yields that
$P\leqst Q$ implies $\int x\,P(\d x)\leq\int x\,Q(\d x)$. Thus
stochastic ordering implies ordering of the expected values but
not vice versa.

There is a vast literature on stochastic orderings, and we only
refer to \cite{MuellerStoyan2002}, \cite{ShakedShanthikumar2007} and \cite{Szekli1995}.

Let $b_{n,p}$ denote the binomial distribution with parameters
$n\in\N$ and $p\in[0,1]$, let $\Poi_\lambda$ denote the Poisson distribution with parameter $\lambda>0$ and let $b_{r,p}^-$ denote the negative binomial distribution (also known as Pascal distribution) with parameters $r\in(0,\infty)$ and $p\in(0,1]$. Recall that
$b_{r,p}^-$ is the probability measure on $\N_0$ with weights
$$b_{r,p}^-(\{k\})={-r\choose k}(-1)^kp^r(1-p)^k={r+k-1\choose k}p^r(1-p)^k\mf k\in\N_0.$$
Further, we denote by
$$\hyp_{B,W,n}(\{k\})={B\choose k}{W\choose n-k}\Big/{B+W\choose n},\qquad k=(n-W)^+,\ldots, B\wedge n$$
the hypergeometric distribution with parameters $B,W\in\N_0$ and $n\in\N$ with $n\leq B+W$.
The main goal of this paper is to prove
necessary and sufficient conditions for $b_{n_1,p_1}\leqst b_{n_2,p_2}$, for $b_{r_1,p_1}^-\leqst b_{r_2,p_2}^-$ and for $\hyp_{B_1,W_1,n_1}\leqst\hyp_{B_2,W_2,n_2}$ in terms of the
parameters $r_1,r_2,n_1,n_2,p_1,p_2,B_1,W_1,B_2,W_2$.

Since stochastic ordering implies ordering of expectations,
$b_{n_1,p_1}\leqst b_{n_2,p_2}$ implies $p_1n_1\leq p_2n_2$, but
this condition is not sufficient for $b_{n_1,p_1}\leqst
b_{n_2,p_2}$. However, if $n:=n_1=n_2$, then
\begin{equation}
\label{E1.1}
b_{n,p_1}\leqst b_{n,p_2} \quad\iff\quad p_1\leq p_2.
\end{equation}
There are various proofs of this statement, the simplest being a
coupling: Let $U_1,\ldots,U_n$ be i.i.d.{} random variables that
are uniformly distributed on $[0,1]$. For $i=1,2$, let
$$N_i=\#\{k:\,U_k\leq p_i\}.$$
Then $N_i\sim b_{n,p_i}$ and $N_1\leq N_2$ almost surely. In Section~\ref{S3} we present a more involved coupling proving
the sufficiency of a characterization of $b_{n_1,p_1}\leqst
b_{n_2,p_2}$ also when $n_1\neq n_2$.

\subsection{The Likelihood Ratio Order}
\label{S1.2}
Before we come to the statement of the main theorem of this article let us briefly discuss a stronger notion of ordering of two probability measures on $\R$, the so-called \emph{monotone likelihood ratio} order.
Let $\mu$ be any $\sigma$-finite measure such that $P$ and $Q$ are absolutely continuous with respect to $\mu$ and $\mu$ is absolutely continuous with respect to $P+Q$. Furthermore, define the respective densities
$$
f=\frac{\d P}{d\mu}\mbu g=\frac{\d Q}{\d\mu}.
$$
$P$ is said to be smaller than or equal to $Q$ in the monotone likelihood ratio order ($P\leqlr Q$) if there exist versions of $f$ and $g$ such that the likelihood ratio
\begin{equation}
\label{Edlr}
x\mapsto \l(x):=\frac{f(x)}{g(x)}\mbsl{is monotone decreasing.}
\end{equation}

Note that the ordering does not depend on the choice of $\mu$; in particular, $\mu=P+Q$ is possible.

It is well known that $P\leqlr Q$ implies $P\leqst Q$ but not vice versa. This will become even more obvious by the following characterization of the monotone likelihood ratio order. Let $\CB(\R)$ denote the Borel $\sigma$-algebra on $\R$. Then we have
\begin{equation}
\label{E1.1b}
P\leqlr Q{\ }\iff{\ }
P(\ARG\given B)\leqst Q(\ARG\given B)\mfa
B\in\CB(\R),\,P(B)>0,Q(B)>0
\end{equation}
by any of \cite[pp.{} 1217-1218]{Pfanzagl1964}, \cite[Theorems 1.1, 1.3]{Whitt1980} or \cite[pp.{} 50-52]{Rueschendorf1991}. In fact, the $\Longleftarrow$ implication is valid even if we replace $\CB(\R)$ by the class of all intervals in $\R$ (see \cite{Pfanzagl1964}) or
by any smaller class $\CC$ of subsets of $\R$ such that for any $r<s$ there exists an $\ve>0$ and a $B\in\CC$ such that $[r-\ve,r]\cup[s,s+\ve]\in B$ (see \cite[Theorem 1.3]{Whitt1980}). In particular, if $P$ and $Q$ live on a discrete subset of $\R$, then it suffices to check the right hand side of \equ{E1.1b} only for sets $B$ of cardinality $2$.

For the binomial distributions, we have $b_{n_1,p_1}\leqlr b_{n_2,p_2}$ if and only if $p_1=0$ or
\begin{equation}
\label{E2.1}
n_1\leq n_2\mbu \frac{n_1p_1}{1-p_1}\leq \frac{n_2p_2}{1-p_2}.
\end{equation}
(See \cite[Theorem~1(iv)]{BolandSinghCukic2002} for a result for a larger class of distributions that comprises the binomial distributions.) In fact, if we exclude the trivial case $p_1=0$, then $n_1\leq n_2$ is clearly necessary for $b_{n_1,p_1}\leqlr b_{n_2,p_2}$. In order to see that \equ{E2.1} is sufficient, assume $n_1\leq n_2$ and let $f_1$ and $f_2$ be the corresponding densities, say with respect to the counting measure on $\N_0$. Then $\l=f_1/f_2$ is decreasing if and only if for all $k=0,\ldots,n_1-1$
$$
1\,\geq\;\frac{f_1(k+1)/f_2(k+1)}{f_1(k)/f_2(k)}
\;=\;\frac{p_1}{1-p_1}\,\frac{1-p_2}{p_2}\,\frac{n_1-k}{n_2-k}.
$$
Clearly, the expression on the right hand side is maximal for $k=0$ and in this case the inequality is equivalent to \equ{E2.1}.

As the monotone likelihood ratio order is stronger than the stochastic order, it is clear that \equ{E2.1} is sufficient for $b_{n_1,p_1}\leqst b_{n_2,p_2}$ but it is not necessary as we will see.

Note that for the Poisson distribution, we have
$$
\Poi_\lambda\leqst\Poi_\mu\quad\iff\quad
\Poi_\lambda\leqlr\Poi_\mu\quad\iff\quad
\lambda\leq\mu.$$
Hence, for this subclass of distributions, stochastic ordering and monotone likelihood ratio ordering coincide.

\subsection{Main Result}
\label{S1.3}
For distributions $P$ and $Q$ on $\N_0$, the likelihood ratio $\l=f/g$ (see \equ{Edlr}) is given by $\l(k):=P(\{k\})/Q(\{k\})$, $k\in\N_0$. Let
\begin{equation}
\label{Edefb}
k_*:=\min(\{k:\,(P+Q)(\{k\})>0\})\mbu k^*=\sup(\{k:\,(P+Q)(\{k\})>0\}).
\end{equation}
If $k^*=\infty$, define $\l(k^*):=\limsup_{k\to \infty}\l(k)$, $\underline\l(k^*):=\liminf_{k\to\infty}\l(k)$ and
the extreme right tail ratio
\begin{equation}
\label{ER1.2c}
\varrho:=\limsup_{k\to \infty}
\frac{P(\{k,k+1,\ldots\})}{Q(\{k,k+1,\ldots\})}.
\end{equation}
If $k^*<\infty$, define $\varrho:=\l(k^*)$. Note that
$\underline\l(k^*)\leq \varrho\leq \l(k^*)$.
In order that $P\leqst Q$ holds, it is clearly necessary that
\begin{equation}
\label{ER1.1}
\l(k_*)\geq 1
\end{equation}
and
\begin{equation}
\label{ER1.2a}
\varrho\leq 1.
\end{equation}
Clearly, \equ{ER1.2a} is implied by
\begin{equation}
\label{ER1.2}
\l(k^*)\leq 1.
\end{equation}
We say that $(P,Q)$ fulfills the left tail condition if \equ{ER1.1} holds and the right tail condition if \equ{ER1.2} holds.

While we have just argued that (at least if $k^*<\infty$ or if $\l(k)$ converges as $k\to\infty$) both tail conditions are necessary for $P\leqst Q$, the next theorem shows that for certain classes of distributions, the tail conditions \equ{ER1.1} and \equ{ER1.2} are in fact equivalent to $P\leqst Q$.
\begin{theorem}
\label{T1}
In each of the following seven cases we have $P_1\leqst P_2$ if and only if the left and right tail conditions hold.

\textbf{(a) Binomial distribution. }
$P_i=b_{n_i,p_i}$ with $p_i\in(0,1)$, $n_i\in\N$, $i=1,2$.\\
Left tail condition:
\begin{equation}
\label{ET1aL} (1-p_1)^{n_1}\geq (1-p_2)^{n_2}.
\end{equation}
Right tail condition:
\begin{equation}
\label{ET1aR} n_1\leq n_2.
\end{equation}

\textbf{(b) Negative binomial distribution. }
$P_i=b^-_{r_i,p_i}$ with $r_i>0$, $p_i\in(0,1]$, $i=1,2$.\\
Left tail condition:
\begin{equation}
\label{ET1bL} p_1^{r_1}\geq p_2^{r_2}.
\end{equation}
Right tail condition:
\begin{equation}
\label{ET1bR} p_1\leq p_2.
\end{equation}
\textbf{(c) Hypergeometric distribution. }
$P_i=\hyp_{B_i,W_i,n_i}$ with $B_1,B_2,W_1,W_2,n_i\in\N_0$, $B_i+W_i\geq n_i\geq1$, $i=1,2$.
Furthermore, assume that
\begin{equation}
\label{ET1c2}
B_2+W_2\,\geq\, B_1+W_1
\end{equation}
or
\begin{equation}
\label{ET1c1}
\{n_1,\,B_1,\,n_2-W_2-1\}\,\cap\,\{n_2,\,B_2,\,n_1-W_1-1\}\neq\emptyset.
\end{equation}
Define $$k_*=(n_1-W_1)^+\wedge(n_2- W_2)^+\mbu k^*=(n_1\wedge B_1)\vee (n_2\wedge B_2).$$
Left tail condition:
\begin{equation}
\label{ET1cL}
\hyp_{B_1,W_1,n_1}(\{k_*\})\geq
\hyp_{B_2,W_2,n_2}(\{k_*\}).
\end{equation}
Right tail condition:
\begin{equation}
\label{ET1cR}
\hyp_{B_1,W_1,n_1}(\{k^*\})\leq
\hyp_{B_2,W_2,n_2}(\{k^*\}).
\end{equation}
\textbf{(d) Hypergeometric versus binomial.} $P_1=\hyp_{B,W,m}$, $P_2=b_{n,p}$ with $B,W,m,n\in\N$, $B+W\geq m$, and $p\in(0,1]$.
Left tail condition:
\begin{equation}
\label{ET1dL}
{W\choose m}\Big/{{B+W\choose m}}\geq(1-p)^{n}.
\end{equation}
Right tail condition:
\begin{equation}
\label{ET1dR}
m\wedge B\,\leq\, n.
\end{equation}
\textbf{(e) Binomial versus hypergeometric. } $P_1=b_{m,p}$, $P_2=\hyp_{B,W,m}$, with $B,W,m\in\N_0$, $B+W\geq m\geq1$, and $p\in[0,1]$.
Right tail condition:
\begin{equation}
\label{ET1eR}
p^m\leq {B\choose m}\Big/{{B+W\choose m}}.
\end{equation}
Left tail condition:
Is implied by the right tail condition.

\textbf{(f) Binomial versus Poisson. } $P_1=b_{n,p}$, $P_2=\Poi_\lambda$ with $n\in\N$, $p\in[0,1]$ and $\lambda>0$.
Left tail condition:
\begin{equation}
\label{ET1fL}
(1-p)^n\geq e^{-\lambda}.
\end{equation}
Right tail condition: Trivially fulfilled.

\textbf{(g) Poisson versus negative binomial. } $P_1=\Poi_\lambda$, $P_2=b^-_{r,p}$ with $p\in(0,1)$ and $r,\lambda>0$.
Left tail condition:
\begin{equation}
\label{ET1gL}
e^{-\lambda}\geq p^r.
\end{equation}
Right tail condition: Trivially fulfilled.
\end{theorem}
For (a), it is obvious that \equ{ET1aL} is the left tail condition and \equ{ET1aR} is right tail condition.

For (b), \equ{ET1bL} is obviously the left tail condition since $b_{r_i,p_i}(\{0\})=p_i^{r_i}$. For the right tail condition, note that for $k\in\N$, we have
$$\left|{-r_i\choose k}\right|=\prod_{l=1}^k\Big(1+\frac{r_i-1}{l}\Big)
\leq\exp\Big(r_i\sum_{l=1}^k\frac1l\Big)\leq e^{r_i}k^{r_i}.$$
Hence
$$\lim_{k\to\infty}\frac{\log(b_{r_i,p_i}^-(\{k,k+1,\ldots\}))}{k}=\log(1-p_i)$$
and the right tail condition \equ{ER1.2a} is equivalent to $p_1\geq p_2$.

For (c) note that $k_*$ and $k^*$ are the minimum and maximum of the support of $\hyp_{B_1,W_1,n}+\hyp_{W_2,B_2,n}$, respectively. Furthermore, note that in the case $n:=n_1=n_2$, condition \equ{ET1c1} is satisfied. In this case the left tail condition simplifies to
\begin{equation}
\label{ET1cL'}
{B_1+W_1-n\choose B_1-k_*}{B_2+W_2\choose B_2}
\geq
{B_2+W_2-n\choose B_2-k_*}{B_1+W_1\choose B_1}
\end{equation}
and the right tail condition becomes
\begin{equation}
\label{ET1cR'}
{B_1+W_1-n\choose B_1-k^*}{B_2+W_2\choose B_2}
\leq
{B_2+W_2-n\choose B_2-k^*}{B_1+W_1\choose B_1}.
\end{equation}

For (d), (e), (f) and (g), the statements are (almost) trivial. In particular, (d) is a consequence of (c) since $b_{n,p}$ is the limit of $\hyp_{\lfloor pN\rfloor,\lfloor(1-p)N\rfloor,n}$ as $N\to\infty$ and for sufficiently large $N$, condition \equ{ET1c2} is satisfied. Taking a further limit we recover (a). Similarly, (e) can be inferred from (c) noting that condition \equ{ET1c1} is satisfied. In Section~\ref{S2} we give the short proofs though, in order to demonstrate the flexibility of our Method 1, described below.

Part (a) of the theorem is not trivial but is not new either. However, in this paper we give new and elementary proofs using different methods.

\textbf{Method 1} is based on likelihood ratio considerations. We show in Proposition~\ref{P2.1.3} that the left and right tail condition are sufficient for stochastic ordering whenever the likelihood ratio $\l$ or $1/\l$ is a unimodal function; that is, if $\l$ is either first monotone increasing and then monotone decreasing or vice versa. In this case we say that $P$ and $Q$ have half-monotone likelihood ratios.

\textbf{Method 2} works for the binomial distribution only and relies on an explicit coupling of two random variables $N_i\sim b_{n_i,p_i}$, $i=1,2$, such that $N_1\leq N_2$ almost surely.

\textbf{Method 3} also works for the binomial distribution only.
Similarly to Method 2, this method is based on the observation that $b_{n,p}$ can be represented as the number of nonempty boxes when we throw a certain random Poisson number of balls into $n$ boxes. Unlike in Method 2,  here we do not construct an explicit coupling of $N_1$ and $N_2$ but give a stochastic comparison of the Markov dynamics of subsequently throwing the balls.

\textbf{Method 4} works for the binomial and negative binomial distribution and relies on explicitly calculating the changes when we modify the parameter $p$ continuously.

\textbf{Method 5} uses infinite divisibility of the negative binomial distribution to give a proof for part (b).
\medskip

\subsection{Organization of the Paper}
In Section~\ref{S1.5} we provide a brief review on stochastic orderings of Bernoulli convolutions. In Sections~\ref{S2} -- \ref{S6}, we give proofs of Theorem~\ref{T1} using the different methods presented above.
\subsection{A Review on Bernoulli Convolutions}
\label{S1.5}
We give a brief review on a result concerning the stochastic ordering of \emph{Bernoulli convolutions} (that comprises part (a) of our Theorem~\ref{T1}) due to Proschan and Sethuraman \cite{ProschanSethuraman1976}.
Fix $n\in\N$ and let
$$
\Delta_n=\big\{\mathbf{p}=(p_1,\ldots,p_n)\in[0,1]^n:\;p_1\geq p_2\geq\ldots\geq p_n\big\}.$$
Let $\mathbf{p}\in\Delta_n$ and let $X_1,\ldots,X_n$ be independent random variables with $\P[X_i=1]=1-\P[X_i=0]=p_i$. Then the distribution of $X_1+\ldots+X_n$ is said to be the Bernoulli convolution $BC_{\mathbf{p}}$ with parameter ${\mathbf{p}}$.

Let ${\mathbf{p}},{\mathbf{q}}\in\Delta_n$. By \cite[Corollary 5.2]{ProschanSethuraman1976} for $BC_{\mathbf{p}}\leqst BC_{\mathbf{q}}$, it is sufficient that
\begin{equation}
\label{E2.2}
\prod_{j=1}^kp_j\;\leq\; \prod_{j=1}^k q_j\mfa k=1,\ldots,n.
\end{equation}
By the obvious symmetry in the problem (changing the roles of the ones and zeros), it is also sufficient to have
\begin{equation}
\label{E2.3}
\prod_{j=k}^n(1-p_j)\;\geq\; \prod_{j=k}^n (1-q_j)\mfa k=1,\ldots,n.
\end{equation}
Note that \equ{E2.2} and \equ{E2.3} are in fact not equivalent.

Assume $n_1,n_2\leq n$ and $p_1=\ldots=p_{n_1}$, $p_{n_1+1}=\ldots=p_n=0$,
$q_1=\ldots=q_{n_2}$, $q_{n_2+1}=\ldots=q_{n}=0$. Then \equ{E2.3} is equivalent to $n_1\leq n_2$ and
$(1-p_1)^{n_1}\geq(1-q_1)^{n_2}$. Hence Theorem~\ref{T1}(a) is a special case of the result of \cite{ProschanSethuraman1976}.

A special case of \cite[Corollary 5.2]{ProschanSethuraman1976} (which is still more general than our Theorem~\ref{T1}(a)) was investigated independently of Proschan and Sethuraman by Ma \cite{Ma1997}. Ma states \cite[Theorem 1]{Ma1997} that if ${\mathbf{q}}\in\Delta_n$ and $p\in(0,1)$, then
\begin{equation}
\label{E2.4}
BC_{\mathbf{q}}\leqst b_{n,p}\quad\iff\quad  b_{n,p}(\{0\})\leq BC_{\mathbf{q}}(\{0\}).
\end{equation}
In fact, the condition on the right hand side of \equ{E2.4} is \equ{E2.3}
(with the roles of $\mathbf{p}=(p,\ldots,p)$ and $\mathbf{q}$ interchanged).
Again, by the obvious symmetry, this statement is equivalent to
\begin{equation}
\label{E2.5}
b_{n,p}\leqst BC_{{\mathbf{q}}}\quad\iff\quad b_{n,p}(\{n\})\leq BC_{\mathbf{q}}(\{n\}).
\end{equation}
Since the hypergeometric distribution is a Bernoulli convolution (see \cite{VatutinMikhajlov1982}), Theorem~~\ref{T1}(d) and (e) could be inferred from \equ{E2.4} and \equ{E2.5}.
A limiting case of \equ{E2.4}, more general than the present
Theorem~\ref{T1}(f), was given in \cite[(A.5)]{Gastwirth1977}.

\Section{Method 1: Half Monotone Likelihood Ratios}
\label{S2}
In this section we provide a criterion which, together with the left and right tail condition (see \equ{ER1.1} and \equ{ER1.2}) is sufficient for stochastic ordering. We first present this method in the general situation and then apply it to all seven cases (a) -- (g) of Theorem~\ref{T1}.
\subsection{A Special Criterion for the Stochastic Order}
\label{S2.1}
\begin{definition}
\label{D2.1.1}
Let $P$, $Q$ be as in Section~\ref{S1.2}. Define the set $\CH$ of pairs $(P,Q)$ such that there exists a version $\l$ of the likelihood ratio $(\d P/\d (P+Q))\big/(\d Q/\d (P+Q))$ with the following properties:
\begin{itemize}
\item[(i)]  There exists an $x_0\in\R$ such that $\l$ is monotone
(increasing or decreasing) on $(-\infty,x_0]$ and is monotone on $[x_0,\infty)$.
\item[(ii)] The left tail and right tail coniditions hold:
\begin{equation}
\label{ED1.1}
\DS\lim_{x\to-\infty}\l(x)\geq 1
\end{equation}
and
\begin{equation}
\label{ED1.2}
\DS\lim_{x\to\infty}\l(x)\leq 1.
\end{equation}
\end{itemize}
If only (i) is fulfilled, then we write $P\hmlr Q$ and say that $P$ and $Q$ have a half monotone likelihood ratio.
\end{definition}

\begin{remark}
\label{R2.1.2}
For distributions $P$ and $Q$ on $\N_0$, the quotient $f/g$ in Definition~\ref{D2.1.1} is the likelihood ratio $\l(k):=P(\{k\})/Q(\{k\})$, $k\in\N_0$. In this case for $P\hmlr Q$ it is sufficient that
\begin{equation}
\label{ER1.3}
\frac{\l(k+1)}{\l(k)}\mbs{is monotone (increasing or decreasing) for }k_*\leq k<k^*.
\end{equation}
That is, \equ{ER1.1}, \equ{ER1.2} and \equ{ER1.3} imply $(P,Q)\in \CH$.
\hfill$\Diamond$\end{remark}
Note that the relation $\hmlr$ is symmetric and reflexive, but it is not transitive.
Furthermore, note that (trivially) $P\leqlr Q$ implies $(P,Q)\in\CH$.

\begin{proposition}
\label{P2.1.3}
If $(P,Q)\in\CH$, then $P\leqst Q$.
\end{proposition}
\textbf{Proof. } For $P=Q$ the statement is trivial. Hence, now assume $P\neq Q$.
Let $x_0$ be as in the definition of $\CH$.
We will show that there exists an $x_1\in\R$ such that $\l(x)\geq 1$ for $x<x_1$ and $\l(x)\leq 1$ for $x>x_1$. Clearly, this implies $P((-\infty,x])\geq Q((-\infty,x])$ for $x<x_1$ and $P((x,\infty))<Q((x,\infty))$ for $x\geq x_1$. Combining these two inequalities, we get $P\leqst Q$.

In order to establish the existence of such an $x_1$, we distinguish three cases.

\emph{Case 1.}\quad If $\l$ is monotone decreasing, then the statement is trivial.

\emph{Case 2.}\quad Assume that $\l$ is monotone decreasing on $(-\infty,x_0]$ and monotone increasing on $[x_0,\infty)$. Hence $\l(x)\geq \l(x_0)$ for all $x\in\R$ which implies $\l(x_0)<1$ unless $\l(x)=1$ for all $x\in\R$ which was ruled out by the assumption $P\neq Q$. By assumption \equ{ED1.2}, we have $\l(x)\leq 1$ for all $x\geq x_0$. Now take $x_1=\sup\{x: \l(x)\geq 1\}\leq x_0$.

\emph{Case 3.}\quad Assume that $\l$ is monotone increasing on $(-\infty,x_0]$. By assumption \equ{ED1.1}, we have $\l(x_0)>1$, $\l(x)\geq 1$ for all $x\leq x_0$ and $\l$ is monotone decreasing on $[x_0,\infty)$. Choose $x_1=\inf\{x\geq x_0:\,\l(x)\leq 1\}$.
\qed
In Sections~\ref{S2.2} -- \ref{S2.4}  we show that any two binomial distributions, negative binomial distributions and hypergeometric distributions with the same sample size parameter, respectively, have half monotone likelihood ratios. For hypergeometric distributions, we can show this also under the assumptions of Theorem~\ref{T1}(c). For other distributions, this method is not applicable in such generality. For example,  $\hyp_{400,509,500}$ and $\hyp_{310,710,700}$ do not have half monotone likelihood ratios. In fact, the likelihood ratio is increasing on $\{0,1,2\}$ (with maximal value $>1$), decreasing on $\{2,\ldots,150\}$ and increasing on $\{150,\ldots,400\}$ (to values $>1$). These distributions are not stochastically ordered as
$$\hyp_{400,509,500}(\{0,\ldots,k\})< \hyp_{310,710,700}(\{0,\ldots,k\})\mf k\leq 44$$
and
$$\hyp_{400,509,500}(\{0,\ldots,k\})> \hyp_{310,710,700}(\{0,\ldots,k\})\mf k\geq 45.$$
On the other hand, it is simple to check numerically that $\hyp_{100,100,18}\leqst\hyp_{21,23,22}$ but that $\hyp_{100,100,18}\not\hmlr\hyp_{21,23,22}$.

It is tempting to try this method also to get a necessary condition for $b_{n,p}\leq \hyp_{B,W,m}$ with $m\neq n$. However, here in general, we do not have $\hyp_{B,W,m}\hmlr b_{n,p}$ as the following example illustrates. Let $\l(k)=\hyp_{21,23,22}(\{k\})/b_{18,0.5106}(\{k\})$. Then $\l(0)=0.0000042$ and $\l$ increases monotonically to $\l(13)=2.05$. Then it decreases to $\l(17)=0.997$ and finally takes the value $\l(18)=1.006$. Although condition (ii) of Definition~\ref{D2.1.1} is fulfilled, we do not have $(\hyp_{21,23,22},b_{18,0.5106})\in\CH$. In fact, $\hyp_{21,23,22}$ and $b_{18,0.5106}$ are not stochastically ordered since
$$\hyp_{21,23,22}(\{0\})-b_{18,0.5106}(\{0\})=-2.5\cdot 10^{-6}<0$$
and
$$
\hyp_{21,23,22}(\{0,\ldots,16\})-b_{18,0.5106}(\{0,\ldots,16\})=8.4\cdot 10^{-8}>0.$$
It is easy to check that $b_{18,1/2}\leqst\hyp_{21,23,22}$ although $\hyp_{21,23,22}\not\hmlr b_{18,1/2}$. In fact, the likelihood quotient $\l(k)$ increases for $k\leq13$, then decreases to $\l(17)=1.393$ and finally takes the value $\l(18)=1.467$.

\subsection{Proof of Theorem~\ref{T1}(a): Binomial Distributions}
\label{S2.2}
\begin{lemma}
\label{L2.2}
Let $n_1,n_2\in\N$ and $p_1,p_2\in[0,1]$. Then $b_{n_1,p_1}\hmlr b_{n_2,p_2}$; that is, $b_{n_1,p_1}$ and $b_{n_2,p_2}$ have half monotone likelihood ratios.
\end{lemma}
\textbf{Proof. }
The cases $p_1\in\{0,1\}$ or $p_2\in\{0,1\}$ are trivial. Hence, now assume $p_1,p_2\in(0,1)$. Furthermore, due to the symmetry of $\hmlr$ we may assume without loss of generality $n_1\leq n_2$.

Denote by
$$\l(k):=\frac{b_{n_1,p_1}(\{k\})}{b_{n_2,p_2}(\{k\})},\quad k=0,\ldots,n_1,$$
the likelihood ratio.
We compute
$$\frac{\l(k+1)}{\l(k)}=\frac{n_1-k}{n_2-k}\,\frac{p_1(1-p_2)}{(1-p_1)p_2}
\mf k=0,\ldots,n_1-1.
$$
Since $n_1\leq n_2$, we see that $k\mapsto \l(k+1)/\l(k)$ is monotone decreasing and hence $\l(k)$ is first monotone increasing and then monotone decreasing.
\qed

\textbf{Proof of Theorem~\ref{T1}(a). }We only have to show sufficiency of the tail conditions \equ{ET1aL} and \equ{ET1aR} for $b_{n_1,p_1}\leqst b_{n_2,p_2}$. By Proposition~\ref{P2.1.3} and Lemma~\ref{L2.2}, it remains to show \equ{ER1.1} and \equ{ER1.2}.
Since we have $n_2\geq n_1$, we have $k_*=0$ and $k^*=n_2$. Since $p_2\geq p_1$, we get $\l(n_2)=b_{n_1,p_1}(\{n_2\})/b_{n_2,p_2}(\{n_2\})\leq 1$; that is, \equ{ER1.2} holds.
Furthermore, by assumption, we have
$$b_{n_1,p_1}(\{0\})=(1-p_1)^{n_1}\geq (1-p_2)^{n_2}=b_{n_2,p_2}(\{0\})$$
which implies \equ{ER1.1}.\qed
\subsection{Proof of Theorem~\ref{T1}(b): Negative Binomial Distributions}
\label{S2.3}

\begin{lemma}
\label{L2.3}
Let $r_1,r_2>0$ and $p_1,p_2\in(0,1]$. Then $b_{r_1,p_1}^-\hmlr b_{r_2,p_2}^-$; that is, $b_{r_1,p_1}^-$ and $b_{r_2,p_2}^-$ have half monotone likelihood ratios.
\end{lemma}
\textbf{Proof. }
The cases $p_1=1$ or $p_2=1$ are trivial. Hence, now assume $p_1,p_2\in(0,1)$. Furthermore, due to the symmetry of $\hmlr$ we may assume without loss of generality $r_1\leq r_2$.

Denote by
$$\l(k):=\frac{b_{r_1,p_1}^-(\{k\})}{b_{r_2,p_2}^-(\{k\})},\quad k\in\N_0,$$
the likelihood ratio.
We compute
$$\frac{\l(k+1)}{\l(k)}=\frac{r_1+k}{r_2+k}\,\frac{1-p_1}{1-p_2}
\mf k\in \N_0.
$$
Since $r_1\leq r_2$, we see that $\l(k+1)/\l(k)$ is monotone increasing. This implies that $\l(k)$ is first monotone decreasing and then monotone increasing; that is $b_{r_1,p_1}^-\hmlr b_{r_2,p_2}^-$.
\qed

\textbf{Proof of Theorem~\ref{T1}(b). }We only have to show sufficiency of the tail conditions \equ{ET1bL} and \equ{ET1bR} for $b^-_{r_1,p_1}\leqst b^-_{r_2,p_2}$. By Proposition~\ref{P2.1.3} and Lemma~\ref{L2.3}, it remains to show show \equ{ER1.1} and \equ{ER1.2}.
Since $p_1\geq p_2$, we get
$$\lim_{k\to\infty}\left(\frac{b_{r_1,p_1}^-(\{k\})}{b_{r_2,p_2}^-(\{k\})}\right)^{1/k}=\frac{1-p_1}{1-p_2}\leq 1$$
which implies \equ{ER1.2}.
Furthermore, by assumption, we have
$$b_{r_1,p_1}^-(\{0\})=(1-p_1)^{r_1}\geq (1-p_2)^{r_2}=b_{r_2,r_2}^-(\{0\})$$
which implies \equ{ER1.1}.\qed

\subsection{Proof of Theorem~\ref{T1}(c): Hypergeometric Distributions}
\label{S2.4}
The left tail condition \equ{ET1cL} implies $(n_1-W_1)^+\leq(n_2-W_2)^+$ and the right tail condition \equ{ET1cR} implies $n_1\wedge B_1\leq n_2\wedge B_2$. Furthermore, trivially we have $P_1\leqst P_2$ and even $P_1\leqlr P_2$ if
\begin{equation}
\label{E3.8.00}
n_1\wedge B_1\leq (n_2-W_2)^+
\end{equation}
(and hence this condition implies the left and right tail condition). Hence in this case, $(P_1,P_2)\in\CH$.
Since this shows the theorem in the case \equ{E3.8.00}, we may henceforth exclude this case. That is, we assume
\begin{equation}
\label{E3.8.01}
(n_2-W_2)^+<n_1\wedge B_1.
\end{equation}

\begin{lemma}
\label{L2.4.1}
Assume that \equ{ET1c1} holds or that \equ{ET1cR} and \equ{ET1c2} hold. Then we have $\hyp_{B_1,W_1,n_1}\hmlr\hyp_{B_2,W_2,n_2}$; that is, $\hyp_{B_1,W_1,n_1}$ and $\hyp_{B_2,W_2,n_2}$ have half-monotone likelihood ratios.
\end{lemma}
\textbf{Proof. }
By the discussion preceeding this lemma, we may assume
\begin{equation}
\label{E3.8.02}
(n_1-W_1)^+\leq (n_2-W_2)^+<n_1\wedge B_1\leq n_2\wedge B_2.
\end{equation}
Now let
$$\begin{aligned}
k_+&:=(n_2-W_2)^+\geq k_*=(n_1-W_1)^+\quad\mbox{and}\\[2mm]
k^+&:=B_1\wedge n_1\hspace{19pt}\leq k^*=B_2\wedge n_2.\end{aligned}
$$
Further, let
$$\l(k):=\frac{\hyp_{W_1,B_1,n_1}(\{k\})}{\hyp_{W_2,B_2,n_2}(\{k\})}$$
with the convention $1/0=\infty$. We have
$$\l(k)\cases{=\infty,&\mfalls k_*\leq k<k_+,\\[2mm]
\in(0,\infty),&\mfalls k_+\leq k\leq k^+,\\[2mm]
=0,&\mfalls k^+<k\leq k^*.
}
$$
For $k\in I:=\{k_*\vee(k_+-1),\ldots,k^*\wedge(k^++1)\}$, we have
$$q(k):=\frac{\l(k+1)}{\l(k)}=\frac{B_1-k}{B_2-k}\;\frac{W_2-n_2+1+k}{W_1-n_1+1+k}\;\frac{n_1-k}{n_2-k}.$$
(Note that $q(k_+-1)=0$ if $k_+> k_*$.)

We are done if we can show that $q(k)-1$ changes the sign at most in $I$. In fact, this implies that, $q(k)$ can cross $1$ at most once.
This in turn implies that $\l$ is half-monotone on $I$ (in the sense of Definition~\ref{D2.1.1}(i)). Since $\l$ is constant on $\{k_*,\ldots,k_+-1\}$ (taking the value $\infty$) and constant on $\{k^++1,\ldots,k^*\}$ (taking the value $0$), we infer that $\l$ is half-monotone on $\{k_*,\ldots,k^*\}$. Hence, by Remark~\ref{R2.1.2}, we get $\hyp_{B_1,W_1,n_1}\hmlr\hyp_{B_2,W_2,n_2}$.

In order to show that $q(k)-1$ changes the sign at most once, we have to rely on the assumption \equ{ET1c2} or \equ{ET1c1}.

Assume first that \equ{ET1c1} holds. There are nine cases to consider and we start with the case $n_1=n_2$. Then for $k\in I$, we have
$$q(k)-1=\frac{B_1(W_2-n_2+1)-B_2(W_1-n_1+1)+[B_1+W_1-B_2-W_2]k}{(B_2-k)(W_1-n_1+1+k)}.
$$
Note that the numerator is affine linear and the denominator is positive for $k\in I$. Hence $q(k)-1$ changes its sign at most once. The other eight cases $n_1=B_2$, $n_2=B_1$, $B_1=B_2$ and so on are similar resulting in an affine numerator and a denominator without sign change.

Now assume that \equ{ET1c2} holds but \equ{ET1c1} does not hold. Then $k^+<k^*$ and
$$q(k)-1=\frac{p(k)}{r(k)}:=\frac{a_2\,k^2\,+\,a_1\,x\,+\,a_0}{(B_2-k)(W_1-n_1+1+k)(n_2-k)}$$
with
$$\begin{aligned}
a_0=&\,B_1n_1W_2-B_1n_1n_2+B_1n_1-B_2n_2W_1+B_2n_2n_1-B_2n_2\\[1mm]
a_1=&\,-B_1W_2+B_1n_2
-B_1+B_1n_1-n_1W_2-n_1\\&\,+B_2W_1-B_2n_1+B_2-B_2n_2+n_2W_1+n_2\\[1mm]
a_2=&\,B_2+W_2-B_1-W_1.
\end{aligned}
$$
For $k\in I$ the denominator is positive. We have
$$q(k^+)-1=-1$$
and hence $p(k^+)<0$. Since $p$ is at most quadratic, condition \equ{ET1c2} (that is, $a_2\geq0$) implies that $p$ changes its sign at most once on $(-\infty,k^+]$. Hence, again $q(x)-1$ changes its sign at most once.
 \qed

\textbf{Proof of Theorem~\ref{T1}(c). }
We only have to show sufficiency of the tail conditions \equ{ET1cL} and \equ{ET1cR} for $\hyp_{B_1,W_1,n}\leqst\hyp_{B_2,W_2,n}$. However, this is an immediate consequence of Proposition~\ref{P2.1.3} and Lemma~\ref{L2.3}.
\qed

\subsection{Proof of Theorem~\ref{T1}(d): Hypergeometric versus Binomial}
\label{S2.5}
For $m\wedge B\,>\,n$, the implications are clear. Hence, without loss of generality, we may and will assume $m\leq n$ and $B\leq n$.

Denoting the likelihood ratio by $\l(k)=\hyp_{B,W,m}(\{k\})/b_{n,p}(\{k\})$, we get that
$$\frac{\l(k+1)}{\l(k)}=\frac{(B-k)(m-k)}{(W-m+k)(n-k)}\,\frac{1-p}{p}\mf k=0,\ldots,(B\wedge m)-1$$
is monotone decreasing and hence $\hyp_{B,W,m}\hmlr b_{n,p}$. It is a simple exercise to check that
$$\hyp_{B,W,m}(\{0\})/b_{n,p}(\{0\})\geq 1\quad \Longrightarrow \quad \hyp_{B,W,m}(\{n\})/b_{n,p}(\{n\})\leq 1$$
and
$$\hyp_{B,W,m}(\{m\})/b_{m,p}(\{m\})\geq 1\quad \Longrightarrow \quad \hyp_{B,W,m}(\{0\})/b_{m,p}(\{0\})\leq 1.$$
Hence the left tail condition \equ{ET1dL} implies $(\hyp_{B,W,m},b_{n,p})\in\CH$ and thus $\hyp_{B,W,m}\leqst b_{n,p}$.
\qed\medskip
\subsection{Proof of Theorem~\ref{T1}(e): Binomial versus hypergeometric}
\label{S2.6}
The proof of Theorem~\ref{T1}(e) is quite similar to the one of part (d). In fact, it is easy to see that the right tail tail condition \equ{ET1eR} implies $(b_{m,p},\hyp_{B,W,m})\in\CH$ and thus $\hyp_{B,W,m}\geqst b_{m,p}$.
\qed\medskip
\subsection{Proof of Theorem~\ref{T1}(f): Binomial versus Poisson}
\label{S2.7}
Clearly, the left tail condition is necessary for $b_{n,p}\leqst\Poi_\lambda$.

Hence now assume that that the left tail condition \equ{ET1fL} holds.
Let
$$\l(k)=\frac{b_{n,p}(\{k\})}{\Poi_\lambda(\{k\})}$$
and compute
$$\frac{\l(k+1)}{\l(k)}=\frac{p}{(1-p)\lambda}\,(n-k)\mf k=0,\ldots,n.$$
Hence $\l(k+1)/\l(k)$ is monotone decreasing and thus $b_{n,p}\hmlr\Poi_{\lambda}$. Since the right tail condition holds trivially and the left tail condition holds by assumption, we infer $(b_{n,p},\Poi_\lambda)\in\CH$ and thus, by Proposition~\ref{P2.1.3}, we get
$b_{n,p}\leqst\Poi_\lambda$.
\qed\medskip

Of course, this result is trivial, since we can even easily derive a coupling:
Let $\hat\lambda=-\log(1-p)\leq\lambda/n$ and let $X_0,X_1,\ldots,X_n$ be independent with $X_i\sim\Poi_{\hat\lambda}$ for $i=1,\ldots,n$ and $X_0\sim\Poi_{\lambda-n\hat\lambda}$. Then
$$S:=X_0+X_1+\ldots+X_n\geq T:=(X_1\wedge 1)+\ldots+(X_n\wedge 1)\mfs$$
and $S\sim\Poi_\lambda$, $T\sim b_{n,p}$.

\subsection{Proof of Theorem~\ref{T1}(g): Poisson versus negative binomial}
\label{S2.8}
Clearly, the left tail condition is necessary for $\Poi_\lambda\leqst b^-_{r,p}$. Furthermore, it is easy to see that the right tail condition always holds.

Hence now assume that that the left tail condition \equ{ET1gL} holds.
Let
$$\l(k)=\frac{\Poi_\lambda(\{k\})}{b^-_{r,p}(\{k\})}$$
and compute
$$\frac{\l(k+1)}{\l(k)}=\frac{\lambda}{(1-p)}\,\frac{1}{k+1}\mf k\in \N_0.$$
Hence $\l(k+1)/\l(k)$ is monotone decreasing and thus $\Poi_{\lambda}\hmlr b^-_{r,p}$. Since the right tail condition holds trivially and the left tail condition holds by assumption, we infer $(\Poi_\lambda,b^-_{r,p})\in\CH$ and thus, by Proposition~\ref{P2.1.3}, we get
$\Poi_\lambda\leqst b^-_{r,p}$.
\qed\medskip

\Section{Method 2: Coupling}
\label{S3}
In this section, we give a proof of Theorem~\ref{T1}(a) that provides an explicit coupling of two random variables $N_i\sim b_{n_i,p_i}$ such that $N_1\leq N_2$ almost surely. Clearly, this implies $b_{n_1,p_1}\leqst b_{n_2,p_2}$.

\textbf{Proof of Theorem~\ref{T1}(a). }
We only have to show sufficiency of the tail conditions \equ{ET1aL} and \equ{ET1aR} for $b_{n_1,p_1}\leqst b_{n_2,p_2}$. Hence, assume  \equ{ET1aL} and \equ{ET1aR}. By
\equ{E1.1}, it suffices to consider the smallest $p_2$ such that
\equ{ET1aL} holds. That is, we may assume
\begin{equation}
\label{E3.1}
(1-p_1)^{n_1}=(1-p_2)^{n_2}.
\end{equation}
Define
$$
\lambda:=-n_1\,\log(1-p_1)\,=\,-n_2\,\log(1-p_2).
$$
For $i=1,2$, let $(X_i(l),l=1,\ldots,n_i)$ be a
family of independent Poisson random variables with parameter $\lambda/n_i$.
(Note that we do not require that $X_1(l_1)$ and $X_2(l_2)$ be
independent.) Then
$$
N_i=\#\big\{l:\,X_i(l)\geq1\big\}\,\sim\,b_{n_i,p_i}.
$$
The idea is to construct a coupling of the $X_i(l)$ such that
\begin{equation}
\label{E3.2}
N_1\leq N_2\mfasts.
\end{equation}
This clearly implies $b_{n_1,p_1}\leqst b_{n_2,p_2}$.

Let $T$ be a Poisson random variable with parameter $\lambda$.
Assume that for $i=1,2$, the family $(F_{i,k},\,k\in\N)$ of random
variables is independent and independent of $T$ and each $F_{i,k}$
is uniformly distributed on $\{1,\ldots,n_i\}$. Then
$$
X_i(l):=\#\big\{k\leq T:\,F_{i,k}=l\big\},\quad l=1,\ldots,n_i,
$$
are independent and Poisson distributed with parameter
$\lambda/n_i$. The remaining task is to construct the families
$(F_{i,k},\,k\in\N)$ such that \equ{E3.2} holds.

For $A_i\subset\{1,\ldots,n_i\}$ let $a_i=\#A_i$ and
$A_i^c=\{1,\ldots,n_i\}\setminus A_i$. For
$r_1\in\{1,\ldots,n_1\}$ and $r_2\in\{1,\ldots,n_2\}$ define
$q^{A_1,A_2}(r_1,r_2)$ depending on whether $a_1<a_2$ or $a_1\geq
a_2$:

If $a_1<a_2$, then let
$$q^{A_1,A_2}(r_1,r_2)=\frac{1}{n_1n_2}.$$
If $a_1\geq a_2$, then let
$$q^{A_1,A_2}(r_1,r_2)=
\cases{
\DS\frac{1}{a_1n_2},&\mfalls r_1\in A_1\mbs{and}r_2\in A_2,\\[4mm]
\DS\frac{a_1n_2-a_2n_1}{a_1n_1n_2(n_2-a_2)},&\mfalls r_1\in
A_1\mbs{and} r_2\in A_2^c,\\[4mm]
\DS\frac{1}{(n_2-a_2)n_1},&\mfalls r_1\in A_1^c\mbs{and}r_2\in A_2^c,\\[3mm]
0,&\msonst. }
$$
Let
$$q^{A_1,A_2}_i(r_i)=\sum_{r_{3-i}=1}^{n_{3-i}}q^{A_1,A_2}(r_1,r_2)$$
denote the $i$-th marginal of $q^{A_1,A_2}$. Clearly, for $a_1<a_2$ we
have $q^{A_1,A_2}_i(r_i)=1/n_i$ for $i=1,2$ and $r_i\in\{1,\ldots,n_i\}$.
Now assume $a_1\geq a_2$. Then for $r_1\in
A_1$,
$$
q^{A_1,A_2}_1(r_1)
=\frac{a_2}{a_1n_2}+(n_2-a_2)\frac{a_1n_2-a_2n_1}{a_1n_1n_2(n_2-a_2)}
=\frac{1}{n_1}.
$$
On the other hand, for $r_1\in A_1^c$,
$$
q^{A_1,A_2}_1(r_1)=(n_2-a_2)\frac{1}{(n_2-a_2)n_1}\,=\,\frac{1}{n_1}.
$$
Analogously, we get for all $r_2\in\{1,\ldots,n_2\}$
$$
q^{A_1,A_2}_2(r_2)=\frac{1}{n_2}.
$$
Thus, independently of the choice of $A_1$ and $A_2$, the
marginals of $q^{A_1,A_2}$ are the uniform distributions on
$\{1,\ldots,n_1\}$ and $\{1,\ldots,n_2\}$, respectively. Now,
define $A_{0,1}=A_{0,2}=\emptyset$. Inductively, choose a pair
$(F_{k,1},F_{k,2})\in\{1,\ldots,n_1\}\times\{1,\ldots,n_2\}$ at
random according to $q^{A_{k-1,1},A_{k-1,2}}$ and define
$A_{k,i}=A_{k-1,i}\cup\{F_{k,i}\}$. Clearly, $a_{T,i}=N_i$, hence
it is enough to show that
\begin{equation}
\label{E3.3}
a_{k,1}\leq a_{k,2}\mfa k\in\N_0.
\end{equation}
For $k=0$, \equ{E3.3} holds trivially. Now we assume that
\equ{E3.3} holds for $k-1$ and we show that it also holds for $k$.
If $a_{k-1,1}<a_{k-1,2}$, then
$$
a_{k,1}\,\leq\, a_{k-1,1}+1\,\leq\, a_{k-1,2}\,\leq\, a_{k,2}.
$$
If $a_{k-1,1}=a_{k-1,2}$, then either $F_{k,1}\in A_{k-1,1}$,
which implies $a_{k,1}=a_{k-1,1}=a_{k-1,2}\leq a_{k,2}$, or
$F_{k,1}\in A_{k-1,1}^c$. In the latter case, according to the
definition of $q^{A_1,A_2}$, we have $F_{k,2}\in A_{k-1,2}^c$,
hence
\[a_{k,1}\,=\,a_{k-1,1}+1\,=\,a_{k-1,2}+1\,=\,a_{k,2}.
\tag*{$\Box$}\]
\Section{Method 3: Markov Chains}
\label{S4}
The aim of this section is to give a proof of Theorem~\ref{T1}(a) that uses the interpretation of the binomial distribution as the distribution of nonempty boxes when we throw successively balls into $n$ boxes. In contrast to Method 2, here we do not construct an explicit coupling of the random variables but use Markov chains in order to get a very quick and elementary proof that could be taught in any first course on probability theory.

Let $n,t\in\N$. Assume that we throw $t$ balls independently into $n$ boxes with numbers $1,\ldots,n$ and denote by $N_{n,t}$ the number of nonempty boxes. Let $T$ be random and Poisson distributed with parameter $\lambda=-n\log(1-p)$. Assume that $T$ is independent of the numbers $N_{n,t}$, $t=1,2,\ldots$.
As indicated in Section~\ref{S3}, the number $N_{n,T}$ is binomially distributed with parameters $n$ and $p$. Hence, in order to show Theorem~\ref{T1}(a), it is enough to show the following proposition.
\begin{proposition}
\label{P4.1}
For each $t\in\N$, the sequence $(N_{n,t})_{n\in\N}$ is stochastically increasing.
\end{proposition}
Proposition~\ref{P4.1} is in fact a special case of a more general result where the probabilities $p_i$ for hitting box $i=1,\ldots,n$ differ from box to box (see \cite{WongYue1973}).
\medskip

\textbf{Proof. }
For each $n$, $(N_{n,t})_{t=0,1,\ldots}$ is a Markov chain on $\{0,\ldots,n\}$ with transition matrix
$$p_n(k,l)=\cases{k/n,&\mfalls l=k,\\[2mm]1-k/n,&\mfalls l=k+1,\\[2mm]0,&\msonst.}$$
Define
$$h_{n,l}(k)=\sum_{j=l}^{n}p_n(k,j)=
\cases{0,&\mfalls k<l-1,\\[2mm]1-k/n,&\mfalls k=l-1,\\[2mm]1,&\mfalls k>l-1,}
$$
and note that $h_{n,l}(k)$ is increasing in $k$ and $n$.

Let $m<n$ and note that trivially $N_{m,0}=0$ is stochastically smaller than $N_{n,0}=0$. By induction, we show that $N_{m,t}\leqst N_{n,t}$ for all $t\in\N_0$. Indeed, for every $\l\in\{0,\ldots,m\}$, by the induction hypothesis and due to the monotonicity of $(n,k)\mapsto h_{n,l}(k)$, we have
$$
\P[N_{m,t+1}\geq l]=\E[h_{m,l}(N_{m,t})]\leq \E[h_{m,l}(N_{n,t})]\leq \E[h_{n,l}(N_{n,t})]=\P[N_{l,t+1}\geq l].
$$
This however implies that $N_{m,t+1}$ is stochastically smaller than $N_{n,t+1}$.
\qed
\Section{Method 4: Analytic Proof}
\label{S5}
The aim of this section is to give proofs of Theorem~\ref{T1}(a) and (b) that rely on changing the parameter $p$ of the distributions continuously and using calculus to compute the dependence of the distributions on this parameter. Although the proofs for (a) and (b) are rather similar, we felt that it is no loss in efficiency to give two separate proofs.

\subsection{Proof of Theorem~\ref{T1}(a): Binomial Distributions}
\label{S5.1}
\setcounter{equation}{0}
We only have to show sufficiency of the tail conditions \equ{ET1aL} and \equ{ET1aR} for $b_{n_1,p_1}\leqst b_{n_2,p_2}$. By
\equ{E1.1}, we only have to consider the case $n_1<n_2$ and $(1-p_1)^{n_1}=(1-p_2)^{n_2}$.

Let $R:=\frac{n_2}{n_1}>1$ and define the map
\begin{eqnarray}              \label{formula for pi(p)}
\pi:[0,1]\to[0,1],\quad p\mapsto 1-(1-p)^{R}.
\end{eqnarray}
Denote by $\pi'(p)=R(1-p)^{R-1}$ the derivative of $\pi$.

For $n\in\N$, $p\in(0,1)$ and $A\subset\{0,\ldots,n\}$ define
$$
b'_{n,p}(A)=\frac{d}{dp}b_{n,p}(A).
$$
Computing the derivative for $A=\{k\}$, $k=0,\ldots,n$, explicitly yields
$$b'_{n,p}(\{k\}) = -n \big[ b_{n-1,p}(\{k\})- b_{n-1,p}(\{k-1\})\big]$$
(where $b_{n-1,p}(\{-1\})=0$ and $b_{0,p}(\{k\})=1$ iff $k=0$). Hence building a telescope sum, we obtain
\begin{equation}
\label{binomial cdf differentiated}
b_{n,p}'(\{0,\ldots,k\})
=
-n\, b_{n-1,p}(\{k\})
\mf n\in \N, \, p\in [0,1],\,k\in \N_0.
\end{equation}

For $k\in\N_0$, define the map
\begin{equation}
f_k:[0,1]\to\R,\quad p\mapsto b_{n_1,\pi(p)}(\{0,\ldots,k\})
          - b_{n_2^{},p}(\{0,\ldots,k\}).
\end{equation}
As $\pi(p_2)=p_1$, we have to show that $f_k(p_2)\geq0$ for all $k$.
Obviously, only the case
$k\in \{1,\ldots,n_1-1\}$ is nontrivial,
and we fix such a $k$ for the rest of this proof.

Since $\pi(0)=0$ and $\pi(1)=1$, we have  $f_k(0)=f_k(1)=0$.
As $f_k$ is differentiable in $(0,1)$ and continuous on $[0,1]$, it is enough to show that
\begin{equation}
\label{ECon1}
f'_k(p)\mbsl{is strictly positive in a neighbourhood of $0$}
\end{equation}
and
\begin{equation}
\label{ECon2} f'_k(p)=0\mbs{for at most one}p\in(0,1).
\end{equation}
Using \eqref{binomial cdf differentiated},  we compute the derivative
\begin{eqnarray*} f_k'(p) &=& \pi'(p) b'_{n_1,\pi(p)}(\{0,\ldots,k\})
          -  b'_{n_2^{},p}(\{0,\ldots,k\})\\
 &=&
n_2^{} b_{n_2^{}-1,p}(\{k\}) - n_1^{}\pi'(p)  b_{n_1^{}-1,\pi(p)}(\{k\})\\
 &=& n_2^{}{n_2^{}-1\choose k}p^k(1-p)^{n_2-1-k}\\
&&\phantom{xxx}    \,-\,n_1R(1-p)^{R-1}{n_1 -1 \choose k}
     \big( 1-(1-p)^R\big)^k(1-p)^{R(n_1-1-k)} \\
 &=&n_2(1-p)^{n_2-1}\cdot
  \bigg[
   {n_2 -1 \choose k}\left(\frac p{1-p}\right)^k
   \,-\,
   {n_1 -1 \choose k}\left(\frac { 1-(1-p)^R}{(1-p)^R}\right)^k
  \bigg].
\end{eqnarray*}
Hence \equ{ECon1} follows from (recall $R=n_2/n_1>1$)
\begin{eqnarray*}                       \label{fprime positive near zero}
 \lim_{p\downarrow 0} \frac {f'_k(p)}{n_2\,p^k}
  &=&   {n_2 -1 \choose k} \,-\,
   {n_1 -1 \choose k}R^k
  \,>\, 0.
\end{eqnarray*}
Now for $p\in\,(0,1)$, we have $f_k'(p) =0$  if and only if
\begin{eqnarray*}                 \label{equation for g}
 g(p) &:=& p(1-p)^{R-1} - a\cdot \big( 1-(1-p)^R\big)
 \,=\, 0
\end{eqnarray*}
where
$$      \label{Def a}
 a:= \left(\frac{{n_1-1 \choose k}}{{n_2-1 \choose k}} \right)^{1/k}.
$$
Since $f'_k(p)>0$ for sufficiently small $p>0$, also $g(p)>0$ for these $p$. Hence in order to show that $g(p)=0$ for at most one $p\in(0,1)$, it is enough to show that $g'(p)=0$ for at most one $p\in(0,1)$. To this end we compute
$$ g'(p) = (1-p)^{R-2} \big[1-aR-(1-a)Rp\big].
$$
Hence $g'(p) =0$ exactly for $p=1$ and $p=\frac{1-aR}{(1-a)R}$. As this shows \equ{ECon2}, the proof is complete.
\qed
\subsection{Proof of Theorem~\ref{T1}(b): Negative Binomial Distributions}
\label{S5.2}
We only have to show sufficiency of the tail conditions \equ{ET1bL} and \equ{ET1bR} for $b^-_{r_1,p_1}\leqst b^-_{r_2,p_2}$. Recall that if $r_1<r_2$ and $X_1$ and $X_2$ are independent random variables with distributions $b^-_{r_1,p}$ and $b^-_{r_2-r_1,p}$, respectively, then $X_1+X_2$ has distribution $b_{r_2,p}$. That is $b^-_{r_1,p}\leqst b^-_{r_2,p}$ if and only if $r_1\leq r_2$. Hence we may assume $r_1\geq r_2$.

\textbf{Step 1. }By a simple computation, we get
\begin{equation}
\label{E4.5}
\frac{d}{dp}b_{r,p}^-(\{0,\ldots,k-1\})=k{-r\choose k}(-1)^k(1-p)^{k-1}p^{r-1}.
\end{equation}
In fact, multiplying both side in \equ{E4.5} by $p^{-r}$, as functions of $r$,
$$p^{-r}\frac{d}{dp}b_{r,p}^-(\{0,\ldots,k-1\})\mbu
k{-r\choose k}(-1)^k(1-p)^{k-1}p^{-1}$$
are polynomials. Hence, it is enough to check \equ{E4.5} for all $r\in\N$. Either by a direct computation or by appealing to the waiting time interpretation that in a Bernoulli chain with success probability $p$, $b_{n,p}^-$ is the distribution of the number of failures before the $n$th success occurs, we get
$$b_{n,p}^-(\{0,\ldots,k-1\})=b_{n+k-1,1-p}(\{0,\ldots,k-1\}).$$
Hence by \equ{binomial cdf differentiated}, we get
$$\frac{d}{dp}b_{n,p}^-(\{0,\ldots,k-1\})=(n+k-1)b_{n+k-2,1-p}(\{k-1\})=k{-n\choose k}(-1)^k(1-p)^{k-1}p^{n-1}$$
as desired.
\smallskip

\textbf{Step 2. }Let $R=r_2/r_1<1$. Fix $k\in\N$. For $p\in(0,1]$ define
$$f(p):=b_{r_1,p^R}^-(\{0,\ldots,k-1\})-b_{r_2,p}^-(\{0,\ldots,k-1\}).$$
It is enough to show that $f(p)>0$ for all $p\in(0,1]$.

Clearly, $f(1)=0$ and $\lim_{p\downarrow0}f(p)=0$. Hence it is enough to show that $f'(p)>0$ for $p$ sufficiently small and $f'(p)=0$ for at most one $p\in(0,1)$. By \equ{E4.5}, we have
$$f'(p)=kp^{r_2-1}\left[{r_1+k-1\choose k}R(1-p^R)^{k-1}-{r_2+k-1\choose k}(1-p)^{k-1}\right].
$$
Hence
$$
\lim_{p\downarrow0}\frac{f'(p)}{kp^{r_2-1}}={r_1+k-1\choose k}R^k-{r_2+k-1\choose k}\,>\,0.
$$
Define
$$g(p):=c\,(1-p^R)-(1-p),$$
where
$$c:=\left(\frac{{r_1+k-1\choose k}R}{{r_2+k-1\choose k}}\right)^{1/(k-1)}.$$
Clearly, $g(p)=0$ if and only if $f'(p)=0$. It is easy to see that
$g(1)=0$ and $g'(p)=0$ if and only if $p=(cR)^{1/(1-R)}$. This implies that $g(p)=0$ for at most one $p\in(0,1)$.
\qed

\Section{Method 5: Infinite Divisibility}
\label{S6}
In this section, for infinitely divisible distributions on $[0,\infty)$, we derive a sufficient criterion (Lemma~\ref{L6.1}) for stochastic ordering in terms of the L{\'e}vy measures. We use this criterion to give a short proof of Theorem~\ref{T1} for negative binomial distributions.

\subsection{Stochastic Ordering of Infinitely Divisibly Laws}
\label{S6.1}
An infinitely divisible distribution $P$ on $[0,\infty)$ is characterized by the deterministic part $\alpha_P\in[0,\infty)$ and the L{\'e}vy measure $\nu_P$ on $(0,\infty)$. The connection is given via the L{\'e}vy-Khinchin formula (see, e.g., \cite[Theorem 16.14]{Klenke2008e})
\begin{equation}
\label{EL-K}
-\log\left(\int e^{-tx}\,P(\d x)\right) =\alpha_Pt+\int \big(1-e^{-tx}\big)\,\nu_P(\d x)\mfa t\geq0.
\end{equation}
If $\mu$ and $\nu$ are two measures on arbitrary measurable spaces, we write $\mu\leq \nu$ if $\mu(A)\leq \nu(A)$ for every measurable set $A$.
If $X$ is a nonnegative infinitely divisible random variables with distribution $\P_X$, we write $\nu_X=\nu_{\P_X}$ and $\alpha_X=\alpha_{\P_X}$ for the corresponding characteristics. If $X$ and $Y$ are independent, then $\nu_{X+Y}=\nu_X+\nu_Y$ and $\alpha_{X+Y}=\alpha_X+\alpha_Y$. Thus for infinitely divisible distributions $P$ and $Q$, we have
\begin{equation}
\label{EEIIST}
\alpha_P\leq \alpha_Q\mbs{and}\nu_P\leq \nu_Q\quad\Longrightarrow\quad P\leqst Q.
\end{equation}
(The opposite implication is not true.)
For example, the negative binomial distribution $b_{r,p}^-$ is infinitely divisible with vanishing deterministic part and with L{\'e}vy measure $\nu_{r,p}$ concentrated on $\N$ and given by
\begin{equation}
\label{Eb-can}
\nu_{r,p}(\{k\})=\lim_{\lambda\downarrow0}\lambda^{-1}b_{\lambda r,p}^-(\{k\})=r\frac{(1-p)^k}{k}\mf k\in\N.
\end{equation}
In fact, a direct computation yields $$-\log\sum_{k=0}^\infty b_{r,p}^-(\{k\})e^{-tk}=r\,\log\big(p^{-1}(1-(1-p)e^{-t})\big)=\sum_{k=1}^\infty \nu_{r,p}(\{k\})\big(1-e^{-tk}\big).$$
Note that $\nu_{r_1,p_1}\leq \nu_{r_2,p_2}$ if $p_1=p_2$ and $r_1\leq r_2$ or if $r_1=r_2$ and $p_1\geq p_2$. Hence, similarly as for the binomial distribution, by \equ{EEIIST}, we get the two relations (for all $p\in(0,1)$ and $r>0$)
\begin{equation}
\label{E4.1}
b_{r_1,p}^-\leqst b_{r_2,p}^-\quad\iff\quad r_1\leq r_2
\end{equation}
and
\begin{equation}
\label{E4.2}
b_{r,p_1}^-\leqst b_{r,p_2}^-\quad\iff\quad p_1\geq p_2.
\end{equation}
For the case where both parameters differ, we need a more subtle criterion.  Let $\mu_1$ and $\mu_2$ be two measures on $(0,\infty)$ with $\mu_i([x,\infty))<\infty$ for all $x\in(0,\infty)$, $i=1,2$. Extending the notion of stochastic ordering to such measures, we write $\mu_1\leqst\mu_2$ if
\begin{equation}
\label{E4.3}
\mu_1([x,\infty))\leq \mu_2([x,\infty))\mfa x\in(0,\infty).
\end{equation}
Note that $\mu_1\leq \mu_2$ implies $\mu_1\leqst\mu_2$.

For a multi-dimensional version of the following lemma, see \cite[Theorem 2.2]{SamorodnitskyTaqqu1993}.
\begin{lemma}
\label{L6.1}
Let $P_i$, $i=1,2$, be infinitely divisible distributions on $[0,\infty)$ with deterministic parts $\alpha_i$ and L{\'e}vy measures $\nu_i$. Assume that
$$\alpha_1\leq \alpha_2\mbu\nu_1\leqst \nu_2.$$
Then there exist random variables $Z_1$ and $Z_2$ with distributions $P_1$ and $P_2$, respectively, such that $Z_1\leq Z_2$ almost surely. In particular, we have $P_1\leqst P_2$.
\end{lemma}
\textbf{Proof. }It is enough to consider the situation $\alpha_1=\alpha_2=0$. Let $$G_i(x)=\nu_i([x,\infty))\mfa x\in(0,\infty)$$
and define the inverse function
$$G_i^{-1}(y):=\inf\big\{x\geq0:\,G_i(x)\leq y\big\}\mf y\in(0,\infty).$$
Let $X$ be a Poisson point process on $(0,\infty)$ with rate $1$. That is, $X$ is an integer valued random measure on $(0,\infty)$, for bounded measurable sets $A$, $X(A)$ is Poisson distributed with the Lebesgue measure of $A$ as parameter, and for pairwise disjoint sets, the values of $X$ are independent random variables. Define
$$Z_i:=\int G_i^{-1}(y)\,X(\d y),\qquad i=1,2.$$
Then for $t\geq0$, we have (compare \cite[Theorem 24.14]{Klenke2008e})
$$\begin{aligned}
-\log\E\big[e^{-tZ_i}\big]=\int_0^\infty \Big(1-e^{-tG_i^{-1}(y)}\Big)\,\d y
&=\int \big(1-e^{-ty})\,\nu_i(\d y)\\
&=-\log\int e^{-tx}\,P_i(\d x).
\end{aligned}
$$
Thus $Z_i$ has distribution $P_i$. By the assumption $G_1\leq G_2$, we have $G_1^{-1}\leq G_2^{-1}$ and hence $Z_1\leq Z_2$ almost surely.
\qed

\subsection{Proof of Theorem~\ref{T1}(b)}
\label{S6.2}
We only have to show sufficiency of the tail conditions \equ{ET1bL} and \equ{ET1bR} for $b^-_{r_1,p_1}\leqst b^-_{r_2,p_2}$.

Recall from \equ{Eb-can} that the negative binomial distribution $b_{r_i,p_i}^-$ is infinitely divisible with deterministic part $\alpha_{r_i,p_i}=0$ and L{\'e}vy measure $\nu_{r_i,p_i}$ being concentrated on $\N$ and given by
$$\nu_{r_i,p_i}(\{k\})=r_i\frac{(1-p_i)^k}{k}\mfa k\in\N.$$
As $p_1\geq p_2$, we have $\nu_{r_1,p_1}\leqlr \nu_{r_2,p_2}$; that is, the map
$$k\mapsto \frac{\nu_{r_1,p_1}(\{k\})}{\nu_{r_2,p_2}(\{k\})}=\frac{r_1}{r_2}\frac{(1-p_1)^k}{(1-p_2)^k}$$
is monotone decreasing. This implies that
$$k\mapsto \phi(k):= \frac{\nu_{r_1,p_1}(\{k,k+1,\ldots\})}{\nu_{r_2,p_2}(\{k,k+1,\ldots\})}
=\frac{r_1}{r_2}\frac{\sum_{l=k}^\infty l^{-1}(1-p_1)^l}{\sum_{l=k}^\infty l^{-1}(1-p_2)^l}$$
is monotone decreasing. By the assumption $p_1^{r_1}\geq p_2^{r_2}$, we have
$$\phi(1)=\frac{r_1\log(p_1)}{r_2\log(p_2)}\leq 1.$$
Hence $\phi(k)\leq 1$ for all $k\in\N$; that is $\nu_{r_1,p_1}\leqst\nu_{r_2,p_2}$ and thus $b_{r_1,p_1}^-\leqst b_{r_2,p_2}^-$ by Lemma~\ref{L6.1}. \qed
\subsection{Proof of Theorem~\ref{T1}(g)}
\label{S6.3}
When viewed from the perspective of infinitely divisible distributions, the statement of Theorem~\ref{T1}(g) is trivial. In fact, $b^-_{r,p}$ is infinitely divisible and the L{\'e}vy measure $\nu_{r,p}$ has total mass $\nu_{r,p}(\N)=-r\log(p)$. Since $\Poi_\lambda$ is infinitely divisible with L{\'e}vy measure $\nu_\lambda=\lambda\delta_1$, we see that $\nu_\lambda\leqst\nu_{r,p}$ if and only if $e^{-\lambda}\geq p^r$.
Hence, the claim follows using Lemma~\ref{L6.1}.
\qed
\section*{Acknowledgement}
We thank Abram M. Kagan and Sergei V. Nagaev
for fruitful discussions which triggered the present work.
\nocite{Pfanzagl1964}
\nocite{Rueschendorf1991}
\def\cprime{$'$}

\end{document}